\documentclass[letterpaper, 10 pt, conference]{ieeeconf}
\IEEEoverridecommandlockouts
\usepackage{amsmath,amsfonts,amssymb}
\usepackage{xcolor}	
\usepackage{graphicx}

\newtheorem{assumption}{Assumption}
\newtheorem{thm}{Theorem}

\newtheorem{rem}{Remark}

\def\R{\mathbb{R}}

\def\N{\mathbb{N}}
\def\K{\mathbf{K}}
\def\Q{\mathbf{Q}}
\def\M{\mathbf{M}}

\def\P{\mathbf{P}}

\def\A{\mathbf{A}}

\def\B{\mathbf{B}}

\def\Y{\mathbf{Y}}
\def\X{\mathbf{X}}

\def\z{\mathbf{z}}

\def\bx{\mathbf{x}}
\def\y{\mathbf{y}}

\def\b{\mathbf{b}}

\def\tx{\tilde{x}}
\def\ttx{\tilde{\mathbf{x}}}

\def\x{\mathbf{x}}

\title{A "joint+marginal" algorithm for polynomial optimization}

\author{Jean B. Lasserre and Tung Phan Thanh
\thanks{J.B. Lasserre is with LAAS-CNRS and the Institute of Mathematics, University of Toulouse, France. {\tt lasserre@laas.fr}}
\thanks{Tung Phan Thanh is with LAAS-CNRS, University of Toulouse, France.
{\tt tphanta@laas.fr} } }

\begin{document}
\maketitle
\begin{abstract}
We present a new algorithm for solving
a polynomial program $\P$ based on the recent "joint + marginal" approach
of the first author for parametric polynomial optimization. The idea is to first consider
the variable $x_1$ as a {\it parameter} and solve
the associated $(n-1)$-variable 
($x_2,\ldots,x_n$) problem $\P(x_1)$
where
the parameter $x_1$ is fixed and takes values in  some interval $\Y_1\subset\R$,
with some probability $\varphi_1$ uniformly distributed on $\Y_1$.
Then one considers the hierarchy of what we call "joint+marginal" semidefinite 
relaxations, whose duals provide a sequence of univariate polynomial 
approximations $x_1\mapsto p_k(x_1)$ that converges to the optimal value function $x_1\mapsto J(x_1)$ of problem $\P(x_1)$, as $k$ increases. Then with $k$ fixed {\it \`a priori}, one computes  $\tx_1^*\in\Y_1$ which minimizes the univariate polynomial $p_k(x_1)$ on the interval $\Y_1$, a convex optimization problem that can be solved via a single semidefinite program.
The quality of the approximation
depends on how large $k$ can be chosen (in general for significant size problems $k=1$ is the only choice).
One iterates the procedure with now an $(n-2)$-variable problem $\P(x_2)$ 
with parameter $x_2$ in some new interval $\Y_2\subset\R$, etc. so as to finally obtain
a vector $\ttx\in\R^n$. Preliminary numerical results are provided.
\end{abstract}

\section{Introduction}~

Consider the general  polynomial program 
\begin{equation}
\label{defpb}
\P:\quad f^*:=\min_\bx\,\{f(\bx)\::\:\bx\in\K\,\}\end{equation}
where $f$ is a polynomial, $\K\subset\R^n$ is a basic semi-algebraic set, and
$f^*$ is the {\it global} minimum of $\P$ (as opposed to a local minimum).
One way to approximate
the global optimum $f^*$ of $\P$ is to solve a hierarchy of either LP-relaxations 
or semidefinite relaxations as proposed in e.g.
Lasserre \cite{lassiopt,lasserrelp}.
Despite practice with the semidefinite relaxations seems to
reveals that convergence is fast, the matrix size in the $i$-th semidefinite relaxation
of the hierarchy  grows up as fast as $O(n^i)$. Hence, for large size
(and sometimes even medium size) problems, only a few relaxations of the hierarchy can be implemented
(the first, second or third relaxation). In that case,
one only obtains a lower bound on $f^*$, and no feasible solution in general.
So an important issue is:

 {\it How can we use the result of the $i$-th semidefinite relaxation to find an approximate 
feasible solution of the original problem?} 

For some well-known special cases of 0/1 optimization like e.g. the celebrated MAXCUT problem, one may generate a feasible solution 
with guaranteed performance, from a randomized rounding procedure
that uses an optimal solution of the first semidefinite relaxation (i.e. with 
$i=1$); see Goemans and Williamson \cite{goemans}. But in general there is no such procedure.

Our contribution is to provide two relatively simple algorithms
for polynomial programs
which builds up
upon the so-called "joint+marginal" approach (in short (J+M)) 
developed in \cite{lasserrepara} for {\it parametric} polynomial
optimization. The (J+M)-approach for variables $\bx\in\R^n$
and parameters $\y$ in a simple set $\Y$,
consists of the standard hierarchy of semidefinite relaxations in \cite{lassiopt}
where one treats the parameters $\y$ also as variables.  But now the moment-approach
implemented in the semidefinite relaxations, considers a {\it joint} probability distribution 
on the pair $(\x,\y)$, with the additional
constraint that the {\it marginal} distribution on $\Y$ is fixed (e.g. the uniform probability distribution on $\Y$); whence the name {\it "joint+marginal"}.

For every $k=1,\ldots,n$, let the compact interval
$\Y_k:=[\underline{x}_k,\overline{x}_k]\subset\R$ be contained in the projection of $\K$ into the $x_k$-coordinate axis.
In the context of the (non-parametric) polynomial optimization
(\ref{defpb}), the above (J+M)-approach can be used as follows in what we call the 
{\bf (J+M)-algorithm}:

$\bullet$ (a) Treat $x_1$ as a parameter in the compact interval
$\Y_1=[\underline{x}_1,\overline{x}_1]$ with associated probability distribution $\varphi_1$ 
uniformly distributed on $\Y_1$.

$\bullet$ (b) with $i\in\N$ fixed, solve the $i$-th semidefinite relaxation of the (J+M)-hierarchy \cite{lasserrepara} applied to problem 
$\P(x_1)$ with $n-1$ variables $(x_2,\ldots,x_n)$ and parameter $x_1$,
which is problem $\P$ with the additional constraint that
the variable $x_1\in\Y_1$ is fixed.
The dual provides a univariate polynomial $x_1\mapsto J^1_i(x_1)$ which,
if $i$ would increase, would converge to $J^1(x_1)$ in the $L_1(\varphi_1)$-norm. (The map $v\mapsto J^1(v)$ denotes the optimal value function of $\P(v)$, i.e. the optimal value of
$\P$ given that the variable $x_1$ is fixed at the value $v$.) 
Next, compute $\tx_1\in\Y_1$, a global minimizer of the univariate polynomial 
$J^1_i$ on $\Y_1$ (e.g. this can be done
by solving a single semidefinite program). Ideally, when 
$i$ is large enough, $\tx_1$ should be close to the first coordinate $x^*_1$of 
a global minimizer $\x^*=x^*_1,\ldots,x^*_n)$ of $\P$.

$\bullet$ (c) go back to step (b) with now $x_2\in\Y_2\subset\R$ instead of $x_1$,
and with $\varphi_2$ being the probability measure uniformly distributed on $\Y_2$. With the same method, compute a global minimizer $\tx_2\in\Y_2$, 
of the univariate polynomial $x_2\mapsto J^2_i(x_2)$ on the interval $\Y_2$.
Again, if $i$ would increase, $J^2_i$ would converge in the $L_1(\varphi_2)$-norm
to the optimal value function $v\mapsto J^2(v)$ 
of $\P(x_2)$ (i.e. the optimal value of $\P$ given that the variable $x_2$ is fixed at the value $v$.) Iterate until one has 
obtained $\tx_n\in\Y_n\subset\R$.

One ends up wih a point $\ttx\in\prod_{k=1}^n\Y_k$ 
and in general $\ttx\not\in\K$.
One may then use $\ttx$ as initial guess of a local optimization
procedure to find a local minimum $\hat{\x}\in\K$. The rational behind the (J+M)-algorithm
is that if $i$ is large enough and $\P$ has a unique global minimizer $\x^*\in\K$,
then $\ttx$ as well as $\hat{\x}$ should be close to $\x^*$.

The computational complexity before the local optimization procedure is less than solving $n$
times the $i$-th semidefinite relaxation in the (J+M)-hierarchy (which is itself 
of same order as the $i$-th semidefinite relaxation in the 
hierarchy defined in \cite{lassiopt}), i.e., a polynomial in the input size of $\P$. 

When the feasible set  $\K$ is convex, one may define the following
variant to obtain a {\it feasible} point $\ttx\in\K$.
Again, let $\Y_1$ be the projection of $\K_1$ into the $x_1$-coordinate axis.
Once $\tx_1\in\Y_1$ is obtained in step (b), consider
the new optimization problem $\P(\tx_1)$ in the $n-1$ variables $(x_2,\ldots,x_n)$,
obtained from $\P$ by fixing the variable $x_1\in\Y_1$ at the value $\tx_1$. 
Its feasible set  is the convex set $\K_1:=\K\cap\{ \x:\,x_1=\tx_1\}$.
Let $\Y_2$ be the projection of $\K_1$ into the $x_2$-coordinate axis.
Then go back to step (b) with now $x_2\in\Y_2$ as parameter and $(x_3,\ldots,x_n)$ as variables, to obtain
a point $\tx_2\in\Y_2$, etc. until a point $\ttx\in\prod_{k=1}^n\Y_k$ is obtained.
Notice that now $\ttx\in\K$ because $\K$ is convex. Then proceed as before with $\ttx$ being the initial guess of a local minimization algorithm to obtain a 
local minimizer $\hat{\x}\in\K$ of $\P$.

\section{The "joint+marginal approach to parametric optimization}

Most of the material of this section is taken from \cite{lasserrepara}.
Let $\R[\bx,\y]$ denote the ring of polynomials in the variables $\bx=(x_1,\ldots,x_n)$, and
the variables $\y=(y_1,\ldots,y_p)$, whereas $\R[\bx,\y]_d$ denotes its subspace of polynomials of degree at most $d$.  
Let $\Sigma[\bx,\y]\subset\R[\bx,\y]$ denote the subset of polynomials that are sums of squares (in short s.o.s.). For 
a real symmetric matrix $\A$ the notation $\A\succeq0$ stands for $\A$ is positive semidefinite.

\subsection*{The parametric optimization problem}

Let $\Y\subset\R^p$ be a compact set, called the {\it parameter} set, and 
let $f,h_j\in\R[\x]$, $j=1,\ldots,m$.
Let $\K\subset\R^n\times\R^p$ be the basic closed semi-algberaic set:
\begin{equation}
\label{set-xy}
\K:=\{(\bx,\y)\,:\, \y\in\Y\,;\: h_j(\bx,\y)\,\geq\,0,\: j=1,\ldots,m\}
\end{equation}
and for each $\y\in\Y$, let 
\begin{equation}
\label{set-x}
\K_\y\,:=\,\{\,\bx\in\R^n\::\: (\x,\y)\,\in\,\K\,\}.
\end{equation}
For each $\y\in\Y$, fixed,  consider the optimization problem: 
\begin{equation}
\label{pb1}
J(\y)\,:=\,\inf_\bx\:\{\,f(\bx,\y)\::\: (\bx,\y)\,\in\,\K\,\}.
\end{equation}
The interpretation is as follows: $\Y$ is a set of parameters and for each instance
$\y\in\Y$ of the parameter, one wishes to compute an optimal {\it decision} vector $\bx^*(\y)$ that solves problem (\ref{pb1}). 
Let $\varphi$ be a Borel probability measure on $\Y$, with a positive density
with respect to the Lebesgue measure on $\R^p$
(or with respect to the counting measure if $\Y$ is discrete). For instance
\[\varphi(B)\,:=\,\left(\int_\Y d\y\,\right)^{-1}\displaystyle\int_{\Y\cap B}d\y,\qquad\forall B\in\mathcal{B}(\R^p),\]
is uniformly distributed on $\Y$.  Sometimes, e.g. in the context of optimization with data uncertainty, $\varphi$ is already specified. 
The idea is to use $\varphi$ (or more precisely, its moments) to 
get information on the distribution of
optimal solutions $\bx^*(\y)$ of $\P_\y$, viewed as random vectors.
In this section we assume that for every $\y\in\Y$, the set $\K_\y$ in (\ref{set-x})  is nonempty.

\subsection{A related infinite-dimensional linear program}

Let $\M(\K)$ be the set of finite Borel probability measures on $\K$, and consider the following 
infinite-dimensional linear program $\P$:
\begin{equation}
\label{pb2}
\rho\,:=\,\inf_{\mu\in\M(\K)}\:\left\{\,\int_\K f\,d\mu \::\: \pi\mu\,=\,\varphi\,\right\},
\end{equation}
where $\pi\mu$ denotes the marginal of $\mu$ on $\R^p$, that is,
$\pi\mu$ is a probability measure on $\R^p$  defined by
$\pi\mu(B):=\mu(\R^n\times B)$ for all $B\in\mathcal{B}(\R^p)$.
Notice that $\mu(\K)=1$ for any feasible solution $\mu$ of $\P$.
Indeed, as $\varphi$ is a probability measure and $\pi\mu=\varphi$ one has
$1=\varphi(\Y)=\mu(\R^n\times\R^p)=\mu(\K)$.

The dual of $\P$ is the the following infinite-dimensional linear program:
\begin{equation}
\label{dual-lp}
\begin{array}{ll}
\rho^*\,:=\,\displaystyle\sup_{p\in\R[\y]}&\displaystyle
\int_\Y p(\y)\,d\varphi(\y)\\
&f(\bx)-p(\y)\,\geq\,0\quad\forall (\bx,\y)\in\K.\end{array}
\end{equation}
Recall that a sequence of measurable functions $(g_n)$ on a measure space 
$(\Y,\mathcal{B}(\Y),\varphi)$ converges to $g$, {\it $\varphi$-almost uniformly},
if and only if for every $\epsilon>0$, there is a set $A\in\mathcal{B}(\Y)$ such that $\varphi(A)<\epsilon$ and
$g_n\to g$, uniformly on $\Y\setminus A$.
\begin{thm}[\cite{lasserrepara}]
\label{th1}
Let both $\Y\subset\R^p$ and $\K$ in (\ref{set-xy}) 
be compact and assume that for every $\y\in\Y$,
the set $\K_\y\subset\R^n$ in (\ref{set-x}) is nonempty.
Let $\P$ be the optimization problem (\ref{pb2}) and let
$\X^*_\y:=\{\bx\in\R^n\,:\,f(\bx,\y)=J(\y)\}$, $\y\in\Y$. Then:

{\rm (a)} $\rho\,=\,\displaystyle\int_\Y J(\y)\,d\varphi(\y)$ and $\P$ has an optimal solution.

\indent
{\rm (b)} Assume that for $\varphi$-almost $\y\in\Y$, the set of minimizers of
$\X^*_\y$ is the singleton 
$\{\bx^*(\y)\}$ for some $\bx^*(\y)\in\K_\y$. Then there is a measurable mapping $g:\Y\to\K_\y$ such that
\begin{equation}
\label{th1-3}
\begin{array}{rcl}
g(\y)&=&\bx^*(\y)\:\mbox{ for every }\:\y\in\Y\\
\rho&=&\displaystyle\int_\Y f(g(\y),\y)\,d\varphi(\y),\end{array}
\end{equation}
and for every $\alpha\in\N^n$, and $\beta\in\N^p$:
\begin{equation}
\label{th1-4}
\int_\K\bx^\alpha\y^\beta\,d\mu^*(\bx,\y)\,=\,\int_\Y\y^\beta \,g(\y)^\alpha\,d\varphi(\y).
\end{equation}

\indent
{\rm (c)} There is no duality gap between (\ref{pb2}) and (\ref{dual-lp}), i.e. $\rho=\rho^*$,
and if $(p_i)_{i\in\N}\subset\R[\y]$ is a maximizing sequence of (\ref{dual-lp}) 
then:
\begin{equation}
\label{cor-dual-1}
\displaystyle \int_\Y\,\vert\,J(\y)-p_i(\y)\,\vert\,d\varphi(\y)\,\to\,0\quad\mbox{as $i\to\infty$}.
\end{equation}
Moreover, define the functions $(\tilde{p}_i)$ as follows:
$\tilde{p}_0:=p_0$, and
\[\y\mapsto\tilde{p}_i(\y)\,:=\,\max\,[\,\tilde{p}_{i-1}(\y),p_i(\y)\,],\quad i=1,2,\ldots\]
Then $\tilde{p}_i\to J(\cdot)$, $\varphi$-almost uniformly.
\end{thm}

An optimal solution $\mu^*$ of $\P$ encodes {\it all} information on the 
optimal solutions $\bx^*(\y)$ of $\P_\y$.
For instance, let $\B$ be a given Borel set of $\R^n$. Then from Theorem \ref{th1},
\[{\rm Prob}\,(\bx^*(\y)\in\B)\,=\,\mu^*(\B\times\R^p)\,=\,
\varphi(g^{-1}(B)),\]
with $g$ as in Theorem \ref{th1}(b).

Moreover from Theorem \ref{th1}(c), any optimal or nearly optimal solution of $\P^*$ provides us with some polynomial lower approximation of the optimal value function $\y\mapsto J(\y)$
that converges to $J(\cdot)$ in the $L_1(\varphi)$ norm. Moreover,
one may also obtain a piecewise polynomial approximation that converges to $J(\cdot)$, $\varphi$-almost uniformly. 

In \cite{lasserrepara} the first author has defined a (J+M)-hierarchy of semidefinite relaxations
$(\Q_i)$ to approximate as closely as desired the optimal value $\rho$. In particular,
the dual of each semidefinite relaxation $\Q_i$ provides a polynomial $q_i\in\R[\y]$
bounded above by $J(\y)$, and $\y\mapsto\tilde{q}_i(\y):= \max_{\ell=1,\ldots i}q_\ell(\y)$ converges
$\varphi$-almost uniformly to the optimal value function $J$, as $i\to\infty$.
This last property is the rationale behind the heuristic developed below.

\section{A "joint+marginal" approach}

Let $\N^n_i:=\{\alpha\in\N^n:\vert\alpha\vert\leq i\}$ with $\vert\alpha\vert=\sum_i\alpha_i$.
With a sequence $\z=(z_{\alpha})$ indexed in the canonical basis
$(\bx^\alpha)$  of $\R[\bx]$, let 
$L_\z:\R[\bx]\to\R$ be the linear mapping:
\[f\:(=\sum_{\alpha}f_{\alpha}(\bx))\,\mapsto\: L_\z(f)\,:=\,
\sum_{\alpha}f_{\alpha}\,z_{\alpha},\qquad f\in\R[\bx].\]

\subsubsection*{Moment matrix}
The moment matrix $\M_i(\z)$ associated with a sequence $\z=(z_{\alpha})$, 
$\alpha\in\N^n_{2i}$, has its rows and columns
indexed in the canonical basis $(\bx^\alpha)$, and with entries.
\[\M_i(\z)(\alpha,\beta)\,=\,L_\z(\bx^{\alpha+\beta})\,=\,z_{\alpha+\beta},\quad
\forall\,\alpha,\beta\in\N^n_i.\]

\subsubsection*{Localizing matrix}
Let $q$ be the polynomial $\bx\mapsto q(\bx):=\sum_{u}q_{u}\bx^u$. The localizing matrix $\M_i(q\,\z)$ associated with 
$q\in\R[\bx]$ and a sequence
$\z=(z_{\alpha})$, has its rows and columns
indexed in the canonical basis $(\bx^\alpha)$, and with entries.
\begin{eqnarray*}
\M_i(q\,\z)(\alpha,\beta)&=&L_\z(q(\bx)\bx^{\alpha+\beta})\\
&=&\sum_{u\in\N^n}q_{u}z_{\alpha+\beta+u},\quad\forall\,\alpha,\beta\in\N^n_i.\end{eqnarray*}
A sequence $\z=(z_{\alpha})\subset\R$ is said to have a {\it representing} finite Borel measure supported on $\K$
if there exists a finite Borel measure $\mu$ such that
\[z_{\alpha}\,=\,\int_\K \bx^\alpha\,d\mu,\qquad\forall\,\alpha\in\N^n.\]

\subsection{A "joint+marginal" approach}

With $\{f,(g_j)_{j=1}^m\}\subset\R[\bx]$,
let $\K\subset\R^n$ be the basic compact semi-algebraic set
\begin{equation}
\label{setk}
\K:=\{\bx\in\R^n\,:\, g_j(\bx)\geq0,\:j=1,\ldots,m\},
\end{equation}
and consider the polynomial optimization problem (\ref{defpb}).

Let $\Y_k\subset\R$ be some interval
$[\underline{x}_k,\overline{x}_k]$, assumed to be contained 
in the orthogonal projection of $\K$ into the $x_k$-ccordinate axis.

For  instance when the $g_j$'s are affine (so that $\K$ is a convex polytope),
$\underline{x}_k$ (resp. $\overline{x}_k$) solves the linear program $\min\:({\rm resp}\,\max\,)\, \{x_k:\,\x\in\K\}$.
Similarly, when $\K$ is convex and defined by concave polynomials, 
one may obtain $\underline{x}_k$ and $\overline{x}_k$,
up to (arbitrary) fixed precision. In many cases, (upper and lower) bound constraints 
on the variables are already part of the problem definition.

Let $\varphi_k$ the probability measure uniformly distributed on $\Y_k$, hence with moments $(\beta_\ell)$ given by:
\begin{equation}
\label{mom1}
\beta_\ell=\int_{\underline{x}_1}^{\overline{x}_1}x^kd\varphi_k(x)\,=\,
\frac{\overline{x}_k^{\ell+1}-\underline{x}_k^{\ell+1}}
{(k+1)(\overline{x}_k-\underline{x}_k)}
\end{equation}
for every $\ell=0,1,\ldots$.
Define the following parametric polynomial program in $n-1$ variables:
\begin{equation}
\label{pb-param1}
J^k(y)\,=\,\min_\bx\,\{f(\bx)\::\: \bx\in\K;\;\:x_k=y\},
\end{equation}
or, equivalently $J^k(y)=\min\,\{f(\bx)\,:\, \bx\in\K_y\}$, where for every $y\in\Y$:
\begin{equation}
\label{setky}
\K_y:=\{\bx\in\,\K;\, x_k=y\}.
\end{equation}
Observe that by definition, $f^*=\displaystyle\min_{x}\{J^k(x):x\in \Y_k\}$, and
$\K_y\neq\emptyset$ whenever $y\in\Y_k$, where $\Y_k$ is the orthogonal 
projection of $\K$ into the $x_k$-coordinate axis.

\subsection*{Semidefinite relaxations}

To compute (or at least approximate)
the optimal value $\rho$ of problem $\P$ in (\ref{pb2}) associated with the parametric optimization problem (\ref{pb-param1}),
we now provide a hierarchy of semidefinite relaxations in the spirit of 
those defined in \cite{lassiopt}. 
Let $v_j:=\lceil({\rm deg\,g_j})/2\rceil$, $j=1,\ldots,m$,
and for $i\geq \max_jv_j$, consider the semidefinite program:
%
\begin{eqnarray}
\label{primal}
\rho_{ik}=&\displaystyle\inf_\z& L_\z(f)\\
\nonumber
&\mbox{s.t.}&\M_i(\z)\succeq0,\:\M_{i-v_j}(g_j\,\z)\succeq0,\quad j=1,\ldots,m\\
\nonumber
&&L_\z(x_k^\ell)=\beta_\ell,\quad \ell=0,1,\ldots 2i,
\end{eqnarray}
where $(\beta_\ell)$ is defined in (\ref{mom1}). We call (\ref{primal}) the {\it parametric semidefinite relaxation} of $\P$ with parameter $y=x_k$.
Observe that without the "moment" constraints $L_\z(x_k^\ell)=\beta_\ell$,
$\ell=1,\ldots 2i$, the semidefinite program (\ref{primal}) is a relaxation of $\P$ and 
if $\K$ is compact, its corresponding optimal value $f^*_i$ converges to $f^*$ as $k\to\infty$; see Lasserre \cite{lassiopt}.

Letting $g_0\equiv 0$, the dual of (\ref{primal}) reads:
\begin{equation}
\label{dual}
\begin{array}{rl}
\rho_{ik}^*=&\displaystyle\sup_{\lambda,(\sigma_j)} 
\sum_{\ell=0}^{2i}\lambda_\ell\,\beta_\ell\\
\mbox{s.t.}&f(\x)-\displaystyle\sum_{\ell=0}^{2i}\lambda_\ell x_k^\ell=
\sigma_0+\displaystyle\sum_{j=1}^m\sigma_j\,g_j\\
&\sigma_j\in\Sigma[\bx],\quad 0\leq j\leq m;\\
&{\rm deg}\,\sigma_jg_j\leq 2i,\quad 0\leq j\leq m.
\end{array}\end{equation}
Equivalently, recall that $\R[x_k]_{2i}$ is
the space of univariate polynomials of degree at most $2i$, and observe that in (\ref{dual}), the criterion reads
\[\sum_{\ell=0}^{2i}\lambda_\ell\,\beta_\ell=\int_{\Y_k}p_i(y)d\varphi_k(y),\]
where $p_i\in\R[x_k]_{2i}$ is the univariate polynomial
$x_k\mapsto p_i(x_k):=\sum_{\ell=0}^{2i}\lambda_\ell x_k^\ell$.
Then equivalently, the above dual may be rewritten as:
\begin{equation}
\label{duall}
\begin{array}{rl}
\rho_{ik}^*=&\displaystyle\sup_{p_i,(\sigma_j)} \int_{\Y_k} p_i d\varphi_k\\
\mbox{s.t.}&f-p_i=\sigma_0+\displaystyle\sum_{j=1}^m\sigma_j\,g_j\\
&p_i\in\R[x_k]_{2i};\:\sigma_j\in\Sigma[\bx],\quad 0\leq j\leq m;\\
&{\rm deg}\,\sigma_jg_j\leq 2i,\quad 0\leq j\leq m.
\end{array}\end{equation}
\begin{assumption}
\label{ass1}
The family of polynomials $(g_j)\subset\R[\x]$ is such that
for some $M>0$, 
\[\x\mapsto M-\Vert\x\Vert^2=\sigma_0+\sum_{j=1}^m\sigma_j\,g_j,\]
for some $M$ and some s.o.s. polynomials $(\sigma_j)\subset\Sigma[\x]$.
\end{assumption}
\begin{thm}
\label{th22}
Let $\K$ be as (\ref{setk}) and Assumption \ref{ass1} hold.
Let the interval $\Y_k\subset\R$ be
the orthognal projection of $\K$ into the 
$x_k$-coordinate axis, and let $\varphi_k$ be the probability measure, uniformly distributed on $\Y_k$. Assume that $\K_y$ in (\ref{setky}) is not empty,
let $y\mapsto J^k(y)$ be as in (\ref{pb-param1}) and
consider the semidefinite
relaxations (\ref{primal})-(\ref{duall}). Then as $i\to\infty$:

{\rm (a)} $\rho_{ik}\uparrow \displaystyle\int_{\Y_k} J^k d\varphi_k$  and $\rho_{ik}^*\uparrow \displaystyle\int_{\Y_k} J^k d\varphi_k$

{\rm (b)} Let $(p_i,(\sigma_j^i))$ be a nearly optimal solution of (\ref{duall}), e.g. such that
$\int_{\Y_k}p_id\varphi_k \geq\rho_{ik}^*-1/i$. Then $p_i(y)\leq J^k(y)$ for all $y\in\Y_k$, and
\begin{equation}
\label{th22-2}
\displaystyle\int_{\Y_k}\vert J^k(y)-p_i(y)\vert\,d\varphi_k(y)\,\to\,0,\quad\mbox{as }i\to\infty.
\end{equation}
Moreover, if one defines $\tilde{p}_0:=p_0$, and 
\[\y\mapsto\tilde{p}_i(y)\,:=\,\max\,[\,\tilde{p}_{i-1}(y),p_i(y)\,],\quad i=1,2,\ldots,\]
then $\tilde{p}_i(y)\uparrow J^k(y)$,  for $\varphi_k$-almost all $y\in\Y_k$, and so 
$\tilde{p}_i\to J^k$, $\varphi_k$-almost uniformly on $\Y_k$.
\end{thm}
\vspace{0.2cm}

Theorem \ref{th22} is a direct consequence of  \cite[Corollary 2.6]{lasserrepara}.

\subsection{A "joint+marginal" algorithm for the general case}

Theorem \ref{th22} provides a rationale for the following (J+M)-algorithm in the general case. In what follows we use the primal and dual semidefinite relaxations
(\ref{primal})-(\ref{dual}) with index $i$ {\it fixed}. 

\vspace{0.2cm}

\noindent
{\bf ALGO 1: (J+M)-algorithm: non convex $\K$, relaxation $i$}
\vspace{0.2cm}

\noindent
{\bf Set $k=1$;}\\
{\bf Step $k$}: {\bf Input:} $\K$, $f$, and the orthogonal projection 
$\Y_k=[\underline{x}_k,\overline{x}_k]$ of $\K$ into the $x_k$-coordinate axis, with associated probability measure $\varphi_k$, uniformly distributed on $\Y_k$.\\
{\bf Ouput:} $\tx_k\in\Y_k$.
\vspace{0.1cm}

\noindent
Solve the semidefinite program (\ref{duall})
and from an optimal (or nearly optimal) solution $(p_i,(\sigma_j))$ of (\ref{duall}), 
get a global minimizer $\tx_k$ of the univariate polynomial $p_{i}$ on 
$\Y_k$.\\
If $k=n$ stop and output $\ttx=(\tx_1,\ldots,\tx_n)$, otherwise set $k=k+1$ and repeat.
\vspace{0.2cm}

Of course, in general the vector $\ttx\in\R^n$ does not belong to
$\K$. Therefore a final step consists of computing a local minimum $\hat{\x}\in\K$, by using some local minimization algorithm
starting with the (unfeasible) initial point $\ttx$. Also note that 
when $\K$ is not convex, the determination of bounds
$\underline{x}_k$ and $\overline{x}_k$ for the interval $\Y_k$ may not be easy, and so one might be forced to use a subinterval $\Y'_k\subseteq\Y_k$ with conservative (but computable) bounds 
$\underline{x}'_k\geq\underline{x}_k$ and $\overline{x}'_k\leq\overline{x}_k$.

\begin{rem}
\label{disconnected}
Theorem \ref{th22} assumes that for every $y\in\Y_k$, the set 
$\K_y$ in (\ref{setky}) is not empty,
which is the case if $\K$ is connected. If $\K_y=\emptyset$ for $y$ in some open subset of $\Y_k$, then
the semidefinite relaxation (\ref{primal}) has no solution ($\rho_{ik}=+\infty$), 
in which case one proceeds by dichotomy on the interval $\Y_k$ until $\rho_{ik}<\infty$.
\end{rem}
\subsection{A "joint+marginal" algorithm when $\K$ is convex}

In this section, we now assume that the feasible set $\K\subset\R^n$ of problem $\P$ is convex (and compact). The idea is to compute $\tx_1$ as in {\bf ALGO 1} and then
repeat the procedure but now for the $(n-1)$-variable problem
$\P(\tx_1)$ which is problem $\P$ in which the variable $x_1$ is {\it fixed} at the value $\tx_1$. This alternative is guaranteed to work if $\K$ is convex (but not always if
$\K$ is not convex).

For every $j\geq 2$, denote by $\bx_j\in\R^{n-j+1}$ the vector $(x_j,\ldots,x_n)$, and by $\ttx_{j-1}\in\R^{j-1}$ the vector $(\tx_1,\ldots,\tx_{j-1})$ (and so 
$\ttx_1=\tx_1$).

Let the interval $\Y_1\subset\R$ be the orthogonal projection of 
$\K$ into the $x_1$-coordinate axis. For every $\tilde{x}_1\in\Y_1$,
let the interval $\Y_2(\ttx_1)\subset\R$ be the orthogonal projection
of the set $\K\cap\{\x :\,x_1=\tilde{x}_1\}$ into the $x_2$-coordinate axis. Similarly,
given $\ttx_2\in\Y_1\times \Y_2(\ttx_1)$,
let the interval $\Y_3(\ttx_2)\subset\R$ be the orthogonal projection
of the set $\K\cap\{\x :\,x_1=\tilde{x}_1;\,x_2=\tilde{x}_2\}$ into the $x_3$-coordinate axis,
and etc. in the obvious way.

For every $k=2,\ldots,n$,  and $\ttx_{k-1}\in\Y_1\times\Y_2(\ttx_1)\cdots \times\Y_{k-1}(\ttx_{k-2})$,
let $\tilde{f}_k(\bx_k):=f((\ttx_{k-1},\bx_k))$, and
$\tilde{g}^k_j(\bx_k):=g_j((\ttx_{k-1},\bx_k))$, $j=1,\ldots,m$. Similarly, let
\begin{eqnarray}
\nonumber
\K_k(\ttx_{k-1})&:=&\{\bx_k\,:\:\tilde{g}^k_j(\bx_k)\geq0,\:j=1,\ldots,m\},\\
\label{kjx}
&=&\{\x_k\,:\: (\ttx_{k-1},\x_k)\in\K\},
\end{eqnarray}
and 
consider the problem:
\begin{equation}
\label{pbpj}
\P(\ttx_{k-1}):\quad \min\,\{\tilde{f}_k(\bx_x)\,:\,\bx_x\in\K_j(\ttx_{k-1})\},\end{equation}
i.e. the original problem $\P$ where the variable $x_\ell$ is fixed at the value
$\tx_\ell$, for every $\ell=1,\ldots,k-1$. 

Write $\Y_j(\ttx_{k-1})=[\underline{x}_k,\overline{x}_k]$, and let $\varphi_k$ be the probability measure uniformly distributed on $\Y_k(\ttx_{k-1})$.

Let $\z$ be a sequence indexed in the monomial basis
of $\R[\x_k]$. With index $i$, fixed, the 
parametric semidefinite relaxation (\ref{primal}) with parameter $x_k$, associated with problem $\P(\ttx_{k-1})$, reads:
\begin{equation}
\label{primalj}
\begin{array}{rl}
\rho_{ik}=\displaystyle\inf_\z& L_\z(\tilde{f}_k)\\
\mbox{s.t.}&\M_i(\z),\:\M_{i-v_j}(\tilde{g}^k_j\,\z)\succeq0,\quad j=1,\ldots,m\\
&L_\z(x_k^\ell)=\beta_\ell,\quad \ell=0,1,\ldots,2i,\end{array}\end{equation}
where $(\beta_\ell)$ is defined in (\ref{mom1}).
Its dual is the semidefinite program (with $\tilde{g}^k_0\equiv 1)$):
\begin{eqnarray}
\label{dualj}
\rho_{ik}^*&=&\displaystyle\sup_{p_{i},(\sigma_j)} 
\int_{\Y_k(\ttx_{k-1})}p_{i}d\varphi_k\\
\nonumber
&\mbox{s.t.}&\tilde{f}_k-p_{i}=
\sigma_0+\sum_{j=1}^m\sigma_j\,\tilde{g}^k_j\\
\nonumber
&&p_{i}\in\R[x_k]_{2i},\:\sigma_j\in\Sigma[\bx_k],\quad j=0,\ldots,m\\
\nonumber
&&{\rm deg}\,\sigma_j\tilde{g}^k_j\leq 2i,\quad j=0,\ldots,m.
\end{eqnarray}
The important difference between (\ref{primal}) and (\ref{primalj}) is the {\it size} of the corresponding semidefinite programs, since 
$\z$ in (\ref{primal}) (resp. in (\ref{primalj})) is indexed in the canonical basis 
of $\R[\x]$ (resp. $\R[\x_k]$).

\subsection*{The (J+M)-algorithm for $\K$ convex}

Recall that the order $i$ of the semidefinite relaxation is fxed. The (J+M)-algorithm consists of $n$ steps.
At step $k$ of the algorithm, the vector $\ttx_{k-1}=(\tx_1,\ldots,\tx_{k-1})$ (already computed) is such that $\tx_1\in\Y_1$ and $\tx_\ell\in\Y_{\ell}(\ttx_{\ell-1})$ for every $\ell=2,\ldots,k-1$,
and so the set $\K_k(\ttx_{k-1})$ is a nonempty compact convex  set.
\vspace{0.2cm}

\noindent
{\bf ALGO 2: (J+M)-algorithm: convex $\K$, relaxation $i$}
\vspace{0.2cm}

\noindent
{\bf Set $k=1$;}\\
{\bf Step $k\geq1$}: {\bf Input:} 
For $k=1$,
$\ttx_0=\emptyset$, $\Y_1(\ttx_0)=\Y_1$; $\P(\ttx_0)=\P$, $f_1=f$ and $\tilde{g}^1_j=g_j$, $j=1,\ldots,m$.\\
For $k\geq2$, $\ttx_{k-1}
\in\Y_1\times\Y_2(\tx_1)\cdots\times\Y_{k-1}(\tx_{k-2})$.\\
 {\bf Output:} $\ttx_{k}=(\ttx_{k-1},\tx_k)$ with $\tx_k\in\Y_k(\ttx_{k-1})$.\\
Consider the parametric semidefinite relaxations (\ref{primalj})-(\ref{dualj}) with parameter
$x_k$, associated with problem $\P(\ttx_{k-1})$ in (\ref{pbpj}).
\begin{itemize}
\item From an optimal solution of (\ref{dualj}), extract the univariate polynomial 
$x_k\mapsto p_{i}(x_k):=\sum_{\ell=0}^{2i}\lambda^*_\ell x_k^\ell$.
\item Get a global minimizer $\tx_k$ of $p_{i}$ on the interval
$\Y_k(\ttx_{k-1})=[\underline{x}_k,\overline{x}_k]$, and set $\ttx_{k}:=(\ttx_{k-1},\tx_k)$.
\end{itemize}
If $k=n$ stop and ouput $\ttx\in\K$, otherwise set $k=k+1$ and repeat.
\vspace{0.2cm}

As $\K$ is convex, $\ttx\in\K$ and one may stop. A refinement is to 
now use $\ttx$ as the initial guess of a local minimization algorithm to
obtain a local minimizer $\hat{\x}\in\K$ of $\P$. In view of Theorem \ref{th22}, the larger the index $i$ of the relaxations (\ref{primalj})-(\ref{dualj}), the better the
values $f(\ttx)$ and $f(\hat{\x})$. 

Of course, {\bf ALGO 2} can also be used when $\K$ is not convex. However, it may happen that at some stage $k$, the semidefinite relaxation (\ref{primalj}) may be infeasible because
$J^k(y)$ is infinite for some values of $y\in\Y_k(\ttx_{k-1})$. This is because the feasible set
$\K(\ttx_{k-1})$ in (\ref{kjx}) may be disconnected.

\section{Computational experiments}

We report on preliminary computational experiments on some
non convex NP-hard optimization problems.
We have tested 
the algorithms on a set of difficult global optimization problems
taken from Floudas et al. \cite{floudas}. To solve the semidefinite programs
involved in {\bf ALGO 1} and in {\bf ALGO 2}, we have 
used the GloptiPoly software \cite{gloptipoly} that implements the hierarchy of semidefinite relaxations defined in \cite[(4.5)]{lassiopt}. 

\subsection{{\bf ALGO 2} for convex set $\K$}
Those problems are taken from \cite[\S 2]{floudas}. The set $\K$ is a convex polytope
and the function $f$ is a nonconvex quadratic polynomial $\x\mapsto \bx' Q\bx+\b'\bx$
for some real symmetric matrix $Q$ and vector $\b$. In Table I
one displays the problem name, the number $n$ of variables, the number $m$ of constraints, the gobal optimum $f^*$, the index $i$ of the semidefinite relaxation in {\bf ALGO 2}, the optimal value obtained using the output
of {\bf ALGO 2} as initial guess in a local minimization algorithm
of the MATLAB toolbox, and the associated relative error.  
As recommended in Gloptipoly \cite{gloptipoly} for numerical stability
and precision, the problem  data have been rescaled to obtain a polytope contained in the box $[-1,1]^n$. As one may see, and excepted for problem 2.8C5,
the relative error is very small.  For the last problem the relative error (about $11\%$) 
is relatively high despite enforcing some extra 
upper and lower bounds $\underline{x_i}\leq x_i\leq \overline{x}_i$,
after reading the optimal solution. However, using $\ttx\in\K$ as initial guess of the local minimization algorithm in MATLAB, one still finds the optimal value $f^*$.
\begin{table}
\label{tab1}
\begin{center}
\begin{tabular}{||r|c|c|r|c|r|r||}
\hline
Prob& $n$& $m$ & $f^*$ & $i$ & ALGO 2 & rel. error\\
\hline
&&&&&&\\
2.2 & 5&11& -17 & 2 &-17.00 & $0\%$\\
2.3 & 6&8& -361.5 & 1 &-361.50 & $0\%$\\
2.6& 10 & 21& -268.01 & 1 &-267.00 & $0.3\%$\\
2.9& 10 & 21& 0 & 1 &0.00 & $0\%$\\
2.8C1& 20 & 30& -394.75 & 1 &-385.30 & $2.4\%$\\
2.8C2& 20 & 30& -884.75 & 1 &-871.52 & $1.5\%$\\
2.8C3& 20 & 30& -8695 & 1 &-8681.7 & $0.15\%$\\
2.8C4& 20 & 30& -754.75 & 1 &-754.08 & $0.09\%$\\
2.8C5& 20 & 30& -4150.41 & 1 &-3678.2 & $11\%$\\
\hline
\end{tabular}
\end{center}
\caption{{\bf ALGO 2} for convex set $\K$}
\end{table}
\addtolength{\textheight}{-3cm} 
\subsection{{\bf ALGO 1} for non convex set $\K$}
Again in Table II below, $n$ (resp. $m$) stands for the number of variables (resp. constraints), and the value displayed in the "{\bf ALGO 1}" column is obtained
in running a local minimization algorithm of the
MATLAB toolbox with the output $\ttx$ of {\bf ALGO 1} as initial guess.

In Problems $3.2$, $3.3$ and $3.4$ from Floudas et al. \cite[\S 3]{floudas},
one has 
$2n$ linear bound constraints and
additional linear and non convex quadratic constraints. As one may see, the results 
displayed in Table II are very good.

For the Haverly Pooling problem 5.2.2 in \cite[\S 5]{floudas}
with three different data sets,
one has $n=9$ and $m=24$ constraints,
among which $3$ nonconvex bilinear constraints and
$18$ linear bound constraints $0\leq x_i\leq 500$,
$i=1,\ldots,9$. In the first run of {\bf ALGO 1} we obtained bad results because the bounds 
are very loose and in the hierarchy of lower bounds $(f^*_k)$ in \cite{lassiopt} 
that converge to $f^*$, if on the one hand $f^*_2=f^*$, on the other hand
the lower bound $f^*_1<f^*$ is loose. In such a case, and in view of the rationale 
behind the "joint+marginal" approach, it is illusory to obtain good results
with {\bf ALGO 1} or {\bf ALGO 2}. Therefore, 
from the optimal solution $\x^*$ in \cite{floudas}, and when $0<x_i^*<500$,
we have  generated stronger bounds $0.4x_i^*\leq x_i\leq 1.6 x_i^*$.
In this case, $f_1^*$ is much closer to $f^*$ and we obtain
the global minimum $f^*$ with {\bf ALGO 1} followed by the local minimization subroutine;
see Table II.  Importantly, in {\bf ALGO 1}, and before running the local optimization subroutine,
one ends up with a non feasible point $\ttx$. Moreover, we had to sometimes use the dichotomy procedure of Remark \ref{disconnected} because if $\Y_k$
is large, one may have $\K_y=\emptyset$ for $y$ in some
open subintervals of $\Y_k$.

Problem 7.2.2 has $13$ linear constraints and $4$ nonlinear constraints with bilinear terms. To handle the non-polynomial function  $x_i^{0.5}$,
one uses the lifting $u_i^2=x_i$, $u_i\geq0$,  $i=5,6$. 
Problem 7.2.6 has only $3$ variables, $6$ linear bound constraints, and
one highly nonlinear constraint (and criterion). Here one uses the lifting
$u^2x_2=1$, $u\geq0$, to handle the term $x_2^{-1}$.
Again one obtains the optimal value $f^*$ with {\bf ALGO 1} followed by a local optimization subroutine.

\begin{table}
\label{tab2}
\begin{center}
\begin{tabular}{||c|c|c|r|c|r|r||}
\hline
Prob& $n$& $m$ & $f^*$ & $i$ & ALGO 1&rel. error\\
\hline
&&&&&&\\
3.2 & 8&22& 7049 & 1 & 7049  &$0\%$\\
3.3 & 5&16& -30665 & 1 & -30665  &$0\%$\\
3.4 & 6&18& -310 & 1 & -298  &$3.8\%$\\
\hline
&&&&&&\\
5.2.2 (1)& 9&24&  400& 1 & 400 &$0\%$\\
5.2.2 (2)& 9&24& 600 & 1 & 600  &$0\%$\\
5.2.3 (3)& 9&24& 750 & 1 & 750 &$0\%$\\
5.2.4& 9&24& 750 & 1 & 750 &$0\%$\\
\hline
&&&&&&\\
7.2.2& 6&17& -0.3746 & 1 & -0.3746 &$0\%$\\
7.2.6& 3&7& -83.254& 1 & -82.3775 &$1\%$\\
\hline
\end{tabular}
\end{center}
\caption{{\bf ALGO 1} for non convex set $\K$}
\end{table}

\subsection{{\bf ALGO 2} for MAXCUT}
Finally we have tested {\bf ALGO 2} on the 
famous NP-hard discrete optimization problem MAXCUT, which consists of minimizing a
quadratic form $\x\mapsto \x' Q\x$ on $\{-1,1\}^n$, for some real symmetric matrix 
$Q\in\R^{n\times n}$. In this case, $\Y_k=\{-1,1\}$ and
the marginal constraint
$L_\z(x_k^\ell)=\gamma_\ell$ in (\ref{primalj}) need only be imposed for $\ell=1$,
because of the constraints $x_k^2=1$ for every $k=1,\ldots,n$.
Accordingly, in an optimal solution of the dual (\ref{dualj}), $p_i\in\R[x_k]$ 
is an affine polynomial $x_k\mapsto p_i(x_k)=\lambda_0+\lambda_1x_k$ for some
scalars $\lambda_0,\lambda_1$. Therefore after solving
(\ref{dualj}) one decides $\tx_k=-1$ if 
$p_i(-1)<p_i(1)$ (i.e. if $\lambda_1>0$) and $\tx_k=1$ otherwise.

Recall that in {\bf ALGO 2} one first compute $\tx_1$, then with $x_1$ fixed at the value
$\tx_1$, one computes $\tx_2$, etc. 
until one finally computes $\tx_n$, and get $\ttx$.
In what we call the "max-gap" variant of {\bf ALGO 2}, one first solves $n$ programs (\ref{primal})-(\ref{dual}) with parameter $x_1$ to obtain an optimal solution $p_i(x_1)=\lambda^1_0+\lambda^1_1x_1$
of the dual (\ref{dual}), then with 
$x_2$ to obtain $(\lambda^2_0,\lambda^2_1)$, etc. finally with
$x_n$ to obtain $(\lambda^n_0,\lambda^n_1)$. 
One then select $k$ such that $\vert\lambda_1^k\vert=\max_\ell\vert\lambda_1^\ell\vert$,
and compute $\tx_k$ accordingly. This is because the larger $\vert\lambda_1\vert$,
(i.e. the larger $\vert p_i(-1)-p_i(1)\vert$), the more likely the choice $-1$ or $1$
is correct.
After $x_k$ is fixed at the value $\tx_k$,
one repeats the procedure for the $(n-1)$-problem $\P(\tx_k)$, etc.

We have tested the "max-gap" variant for
MAXCUT problems on random graphs with $n=20, 30$ and $40$ nodes.
For each value of $n$, we have solved $50$ randomly generated
problems and $100$ for $n=40$. The probability
$\varphi_k$ on $\Y_k=\{-1,1\}$ is uniform (i.e.,
$\beta_1=0$ in (\ref{primalj})).
Let $f^*_1$ denote the optimal value of 
the Shor's relaxation with 
famous Goemans and Williamson's 0.878 performance guarantee.
Let $\rho$ denote the cost of the solution $\x\in\{-1,1\}^n$ generated by
the {\bf ALGO 2}. In Table III we have reported the average relative error
$(\rho-f^*_1)/\vert f^*_1\vert$, which 
as one may see, is comparable
with the Goemans and Williamson (GW) ratio. 

\begin{table}
\label{tab3}
\begin{center}
\begin{tabular}{|c|c|c|c|}
\hline
n& 20 & 30 & 40\\
\hline
&&&\\
$(\rho-f^*_1)/\vert f^*_1\vert$&10.3\% & 12.3\%&12.5\% \\
&&&\\
\hline
\end{tabular}
\end{center}
\caption{Relative error for MAXCUT}
\end{table}

\section{Conclusion}

First preliminary results are promising, even with small relaxation order $i$. When the 
feasible set is non convex, it may become difficult to obtain a feasible
solution and an interesting issue for further investigation is
how to proceed when $\K_y=\emptyset$ for $y$ in some
open subinterval of $\Y_k$ (proceeding by dichotomy on $\Y_k$ is one possiblity).

\end{document}